\documentclass[12pt,twoside,leqno]{amsart}  

\usepackage{amsmath,amscd,amssymb,amsfonts}

\newcommand{\sF}{{\mathcal F}}
\newcommand{\sG}{{\mathcal G}}
\newcommand{\sg}{{\it g}}
\newcommand{\sH}{{\mathcal H}}

\newcommand{\sL}{{\mathcal L}}

\newcommand{\sN}{{\mathcal N}}
\newcommand{\sO}{{\mathcal O}}

\newcommand{\st}{{\it t}}

\newcommand{\C}{{\mathbb C}}

\newcommand{\G}{{\mathbb G}}
\renewcommand{\H}{{\mathbb H}}

\newcommand{\pP}{{\mathbb P}}

\newcommand{\Z}{{\mathbb Z}}

\newcommand{\p}{{\bf P} }
\newcommand{\z}{{\bf G/ P} }
\newcommand{\Tt}{\tau }
\newcommand{\Dif}{\Omega }
\newcommand{\dif}{\omega }

\newcommand{\lra}{\longrightarrow }
\newcommand{\Lra}{\Longrightarrow }

\newcommand{\Hom}{{\rm Hom}}

\newcommand{\Spec}{{\rm Spec \,}}

\theoremstyle{plain}
\newtheorem{thm}{Theorem}
\newtheorem{lem}[thm]{Lemma}

\newtheorem{prop}[thm]{Proposition}
\newtheorem{remark}[thm]{Remark}

\numberwithin{thm}{section}
\numberwithin{equation}{section}
\def\pf{\noindent{\sc Proof. }}
\def\qed{\rm Q.E.D.}

\begin{document}

\title{ Infinitesimal  Torelli theorem for cyclic coverings of generalized flag varieties I}

\author{Pedro Luis del Angel}
\address{Pedro Luis del Angel} 
\address{\vskip-6.5mm Cimat, Guanajuato, Mexico}
\email{luis@cimat.mx}

\author{Herbert Kanarek }
\address{Herbert Kanarek }
\address{\vskip-6.5mm Famat, U. de Guanajuato, Guanajuato, Mexico}
\email{herbert@cimat.mx}

\date{}

\begin{abstract} We give an effective  infinitesimal Torelli theorem 
for cyclic  covers  of  $\z$ where ${\bf G}$ is a simple algebraic group and ${\bf P}$ is a maximal parabolic subgroup.
\end{abstract}
\thanks{Partially supported  by CONCYTEG Grant 04-02-K121-024 (2), CONACyT Grant P083512 and by the University Duisburg-Essen}

\maketitle

\section*{Introduction}

The Torelli problem for a given family of varieties ask wheter 
varieties of the family can be distinguished by their Hodge 
structures. It is known to be the case for curves \cite{torelli},
\cite{andreoti} and \cite{weil}, some K3 surfaces \cite{shapiro}
and Prym varieties \cite{kanev}. 
\\ \\
There are some variants that have been studied: whether the period  
map is an immersion (local Torelli), wheter its differential 
is injective in the deformation space (infinitesimal Torelli) and 
whether the  map is  
generically injective (generic Torelli). Generic Torelli 
is known to hold 
for hypersurfaces of degree $d$ in the projective space (except for 
few exceptions, see \cite{green}), while the infinitesimal Torelli
is known for hypersurfaces of high degree of $\z$ \cite{plar2}.
\\
Let $X$ be a generalized flag variety $\z$, with ${\bf G}$ a simple algebraic group and ${\bf P}$ a maximal parabolic subgroup, embedded 
minimally and equivariantly into a projective space by an ample line bundle 
$\sO_X(1)$ (which generates the Picard group of $X$ in this case)
and let $\Omega ^q(k) = \Omega ^q \otimes \sO_X(1)^k$ be the sheaf of holomorphic
$q-$forms
on $X$ tensored with the $k^{th}$ power of $\sO_X(1)$.
\ \\ \ \\
Let $Z$ be a simple covering of $X$ of degree $N$, in particular $Z$ can be 
embeded in 
$\Spec_{\sO_X}(\oplus_{i=0}^{\infty}\sL^{-i})$ for some ample line bundle 
$\sL=\sO(k)$.
Following some ideas of Flenner \cite{flenner}, Green \cite{green} and 
Ivinskis 
\cite{ivinskis},
we want to show that the infinitesimal Torelli theorem holds for $Z$ under 
certain effective
conditions (see \ref{torelli} ) on the degree of the 
covering and the degree of $\sL$.
\ \\ \ \\
The infinitesimal Torelli theorem ask about the
injectivity of the tangent map
\begin{eqnarray*}
T:H^1(Z,\tau_Z)\to T_{m(0)}D=\oplus_p 
\Hom(H^{n-p}(Z,\Dif^p_Z),H^{n-p+1}(Z,\Dif^{p-1}_Z))
\end{eqnarray*}
(see below), but this can be reduced (as in \cite{ivinskis}) to 
the vanishing of some cohomology groups on the base space $X$, which can be
done in our case using an idea of Deligne-Dimca \cite{dd}, the results of Snow 
\cite{snow1}and \cite{snow2} as well as a classical theorem of Bott 
\cite{bott}.
\\ \\
More precisely, let $B$ be the Kuranishi space of deformations of $Z$. 
Suppose B is smooth 
and let $m:B\to D$ be the local period map into the period domain $D$
of all Hodge structures on $H^k(Z,\C)$, where $k=\dim Z$. The map $m$ is
holomorphic, with tangent map at the point $Z$:
$$T:H^1(Z,\tau_Z)\to T_{m(0)}D$$
We say that the infinitesimal Torelli theorem holds for $X$ if the 
tangent map $T$ is injective. 
\\ \\
The paper is organized as follows: In section 1 we remember some 
basic facts about generalized flag varieties with Picard group $\Z$, in particular 
its  classification. In section 2 we recall some vanishing theorems
for the varieties corresponding to Lie  Groups 
of  types $A, B, C, D, E_{6}$ and $ E_{7}$. 
Section 3 reviews the definition of a simple
cover and some basic facts about line bundles on them. In section 4 we
prove  the Infinitesimal  Torelli theorem for cyclic covers of $\z$ under the assumptions
above (i.e., $Pic(\z)\cong\Z$).
\ \\
\section{Generalized flag varieties having Picard group $\Z$}

Any generalized flag variety is of the form $X:=\z$ with
${\bf G}$ a simple algebraic group and ${\bf P}$ a parabolic proper 
subgroup. If ${\bf T}$ is a maximal torus of ${\bf G}$, one associates to it a set $\Phi$ of characters of the torus, called  the {\it root system} of ${\bf G}$ relative to ${\bf T}$ (or just the root system of ${\bf G}$), together with a base $\Delta\subset\Phi$ of the vector space generated by the characters of the torus. 
\\ \\
Since $\Delta=\{\alpha_1, \cdots , \alpha_l\}$ is a base of the vector space generated by the characters of ${\bf T}$, then any root $\alpha\in\Phi$ can be writen as $\alpha=\sum_i a_i\alpha_i$. In particular a root is said to be positive if  $a_i\ge 0$ for all $i$. The set of positive roots is denoted as $\Phi^+$ and one has $\Phi=\Phi^+ \cup \Phi^-$, where $\Phi^+ = -\Phi^-$.
\\ \\
The tangent space to ${\bf G}$ at the identity is a Lie algebra $\sg$ and the tangent space to
$\p$ at the identity is a Lie subalgebra ${\it p}$. If ${\bf T}$ is a maximal torus on ${\bf G}$ contained in $\p$, then 
\\ \\
\begin{eqnarray*}
\sg = \st \oplus\left({\displaystyle \bigoplus_{\alpha\in\Phi^-}}\ \sg_{\alpha}\right)\oplus\left({\displaystyle \bigoplus_{\alpha\in\Phi^+}}\ \sg_{\alpha}\right)
\end{eqnarray*}
and 
\begin{eqnarray*}
{\it p} = \st \oplus\left({\displaystyle \bigoplus_{\alpha\in\Phi^-}}\ \sg_{\alpha}\right)\oplus\left({\displaystyle \bigoplus_{\alpha\in\Phi_c}}\ \sg_{\alpha}\right)
\end{eqnarray*}
\\ \\
where  $\st$ is the Lie algebra associated to the torus $T$, $\sg_{\alpha}$ is a one dimensional vector space associated to the root $\alpha$ and $\Phi_c$ is a subset of $\Phi^+$ which depends on $\p$, called the set of {\it positive compact roots}. In particular it implies that the tangent space of $\z$ at any point is isomorphic to
\\
\begin{eqnarray*}
\Tt = {\displaystyle \bigoplus_{\alpha\in\Phi^+-\Phi_c}}\ \sg_{\alpha}
\end{eqnarray*}
\\ 
The parabolic group $\p$ is said to be generated by $A\subset\Delta=\{\alpha_1, \cdots , \alpha_l\}$ if 
\begin{eqnarray*}
\Phi_c =\{\alpha=\sum_i a_i\alpha_i \in\Phi^+\; |\; a_j=0 \;\mbox{for all}\; \alpha_j\in A\}
\end{eqnarray*}
Any parabolic group is generated by some $A\subset\Delta$. If $A$ consist of a single element, then we say that $P$ is a {\it maximal parabolic subgroup} generated by that element. Observe that 
\begin{eqnarray*}
\Tt = {\displaystyle \bigoplus_{\alpha\in\Phi^+-\Phi_c}}\ \sg_{\alpha}
\end{eqnarray*}
where $\Phi^+-\Phi_c$ consists of the positive roots for which the coefficient of $\alpha_j$ is strictly positive, for every
$\alpha_j\in A$. 
\\ \\
The simple algebraic groups are completely
classified and they have type: $A_l,  B_l, C_l, D_l, E_6$, $E_7$, $E_8$, $F_4$ and $G_2$. We will only consider maximal parabolic subgroups and the clasical groups ($A_l, B_l, C_l$ and $D_l$ ) in this paper, but the results are also true for  the special groups as well.
\\ \\
1.1) If the type of ${\bf G}$ is $A_l$, then ${\bf G} = SL(l+1)$, 
${\bf P}$ is generated by
$\{ \alpha _r \}$ for some $1\le r \le l+1 = r + t$ and $X=\G r(r,t)$.
\\ \\
1.2) If the type of If the type of ${\bf G}$ is $B_l$, $l\ge 2$, then ${\bf G}=SO(2l+1)$.
\\ \\
1.3) If the type of ${\bf G}$ is $C_l$, then ${\bf G}=Sp(2l)$ .
\\ \\
1.4) If the type of ${\bf G}$ is $D_l$, then ${\bf G}=SO(2l)$. 
\\  \\
The following is a well known fact on the theory of algebraic groups, see for instance (\cite{borel}, section 11) or (\cite{humphreys}, section 16).
\\ \\
\begin{lem}
Let  ${\bf G}, \p$ and  ${\bf T}$ be as before, and let $N_{{\bf G}}({\bf T})$ and $C_{{\bf G}}({\bf T})$ be the normalizer and the centralizer of ${\bf T}$ on ${\bf G}$ respectively. Then the group
\begin{eqnarray*}
{\displaystyle W:= \frac{N_{{\bf G}}({\bf T})}{C_{{\bf G}}({\bf T})} }
\end{eqnarray*}
is finite. It is called the {\it Weyl group} of ${\bf G}$ (relative to ${\bf T}$).
\end{lem}
\ \\
\begin{remark}
There is a well defined notion of {\it length} for the elements of the Weyl Group, which we do not want to make precise here because it will not be really necessary for what comes next. Let us just mention that, in the case of groups of type $A_l$, the Weyl group is the group of permutations $S_l$ and, in this case, the length of an element is the usual length of a permutation.
\end{remark}
\ \\  \ \\
\section{The cohomology of $\Omega _X^q(k)$}

Associated to every root $\alpha_i\in\Delta$ there exist a fundamental weight 
$\lambda_i$ wich satisfies 
${\displaystyle \frac{2(\alpha_i, \lambda_j)}{(\alpha_i,\alpha_j)} = \delta_{i,j}}$ for
every $\alpha, \beta \in \Delta$, where $( , )$ is the killing form on the span of $\Delta$. 
\\ \\
Let $\delta = \sum_{\alpha\in\Delta} \lambda_{\alpha}$ and for every $\beta \in\Phi$ let
$ht (\beta):=(\beta, \delta)$ be the height of $\beta$.
\\ \\
Given a weight $\lambda$, we say that $\lambda$ is {\it singular} if
there exist a positive root $\alpha$ such that $(\alpha, \lambda)=0$. If
$\lambda$ is not singular, we say that it is {\it regular of index $p$},
where $p$ is the number of positive roots $\alpha$ for wich 
$(\alpha, \lambda)<0$.

\begin{thm}
(Bott \cite{bott}, theorems {\rm iv} and {\rm iv'}) Let $\sF$ be an homogeneous vector bundle over $\z$, which is defined by an irreducible representation of ${\bf P}$ with highest weight $\lambda$.
\\
1) If $\lambda + \delta$ is singular, then
$H^k(\z, \sF)=0$ for all $k$, 
\\
2) If $\sF$ is regular of index $p$, then $H^k(\z,\sF)=0$ for every $k\ne p$.
\end{thm}
 
\hfill $\Diamond $

\begin{remark}
If the maximal parabolic subgroup ${\bf P}$ is generated by $\{\alpha_j\}$, then the ample line bundle $\sL $ which generates $Pic(\z)$ is defined by an irreducible representation of ${\bf P}$ with highest weight $\lambda_j$.
\end{remark}

Let $\Phi_c$ is the set of positive compact roots as before and let
\begin{eqnarray*}
W_1=\{\omega\in\; W\; |\; \omega^{-1}\alpha > 0 \; \mbox{for all}\; \alpha\in\Phi_c\},
\end{eqnarray*}
 where $W$ is the Weyl group of ${\bf G}$. Define
\begin{eqnarray*}
W_1(q):=\{\omega\in W_1\; |\; length(\omega)=q\}.
\end{eqnarray*}
The highest weights of the fully reducible ${\bf P}$-module $\Omega^q_{\z}(k)$ occur with multiplicity one and are precisely the weights of the form
\begin{eqnarray*}
\omega\delta-\delta+k\lambda_j, \hskip2cm{} \omega\in W_1(q),
\end{eqnarray*}
so, in order to determine which cohomology groups $H^i(\z,\Omega^q(k))$ vanish, we only need to check whether $\omega\delta+k\lambda_j$ is singular or not, and to compute its index of regularity when it is not. \\ \\
We can further simplify the problem observing that for every positive compact root $\alpha$ one has
$(\alpha,\lambda_j)=0$, therefore
\begin{eqnarray*}
(\omega\delta+k\lambda_i, \alpha) =(\omega\delta,\alpha)=(\delta,\omega^{-1}\alpha)>0
\end{eqnarray*}
since $\alpha\in\Phi_c$ and $\omega\in W_1$. \\ \\
Now, if $\alpha$ is a positive non-compact root, then $\alpha$ involves $\alpha_j$ with coefficient $1$
and $(\lambda_j,\alpha)=(\lambda_j,\alpha_j)=(\alpha_j,\alpha_j)/2 =c$, thus 
$\omega\delta+k\lambda_j$ is singular if and only if there is an $\alpha\in\Phi^+ - \Phi_c$ such that
$c\; k=-(\delta,\omega^{-1}\alpha)$. \\ \\
If $\omega\delta+k\lambda_j$ is not singular, then its index of regularity is
\begin{eqnarray*}
p=| \;\{\alpha\in\Phi^{+}-\Phi_c\; |\; ck<-(\delta,\omega^{-1}\alpha )\}\; |.
\end{eqnarray*}
\begin{remark}
$c=1$, except when ${\bf G}$ is of type $BI$ or $CI$, where $c=2$.
\end{remark}
\begin{remark} \label{canonico}
With the above notation, since $\omega_{\z}$ is a line bundle, it is an irreducible $P$-module associated
to the fundamental weight ${\displaystyle \delta-\sum_{\alpha\in\Phi_c}\alpha}=-d_0\lambda_j$.
\end{remark}
It is not hard to prove that
\begin{lem} \label{bott1}
If ${\displaystyle \mu=\hbox{max}\;\{\frac{(\alpha,\delta)}{c}\; |\; \alpha\in \Phi^+\}}$ and
$k>\mu$ or $k>q$, then 
\begin{eqnarray*}
H^i(\z,\Omega^q(k))=0 
\end{eqnarray*}
for all $ i>0$.
\end{lem}

\section{Simple coverings}

Let $\sL$ be an ample divisor on $X$ and consider the varieties  
$S:=\Spec_{\sO_X}(\oplus_{i=0}^{\infty}\sL^{-i})$ and \label{coverings}
$\bar{S}=\pP(\sO_X\oplus\sL^{-1})$. We say that a smooth variety $Z$ is a 
simple 
covering of $X$ of degree $N$ if there is a finite map $f:Z\to X$ of degree 
$N$ and 
an embedding $Z\stackrel{i}{\to} S\hookrightarrow \bar{S}$ for some ample line 
bundle $\sL$ which makes the following diagram commute:
\begin{eqnarray*}
\begin{matrix}
Z & \stackrel{i}{\to} & S                     & \hookrightarrow           & \bar{S} \cr
  & {\small f}\searrow & {\small \pi}\downarrow & \swarrow {\small\bar{\pi}} &  \cr
  &                   & X                     &                           &  
\end{matrix}
\end{eqnarray*}
where $\pi$ and $\bar{\pi}$ are the natural projections onto $X$.
\ \\
\begin{lem}(\cite{ivinskis}, lemma 1.2) \label{ivinskis1} If $f:Z\lra X$ is a
simple covering of degree $N$, then \\
1) $f_*\sO_Z=\oplus_{i=0}^{N-1}\sL^{-i}$,\ \\
2) $\omega_{Z/X}= f^*\sL^{N-1}$.
\end{lem}

\begin{lem}(\cite{ivinskis}, lemma 2.2) \label{ivinskis2}
Let $Z$ be an smooth projective variety.  The following are equivalent:\ \\
1) $Z$ is a simple covering of $X$ of degree $N$ with respect to the 
ample line bundle $\sL$. \ \\
2) $Z$ is isomorphic to the zero set of a section in $H^0(\bar{S}, 
\sO_{\bar{S}}(N)\otimes\bar{\pi}^*\sL^N).$
\end{lem}
\ \\ \\
Adjunction formula gives us also $\omega_Z= f^*(\sL^{N-1}\otimes\omega_X)= f^*(\sL^{N-1}\otimes \sO_X(-d_0))$, in view of remark \ref{canonico}.
\ \\ \ \\
The following is a well know result which we state in our particular situation
\ \\
\begin{lem}(\cite{hartshorne}, Ex. III, 8.4) \label{projectivebundle}
\begin{eqnarray*}
R^i\bar{\pi}_*\sO_{\bar{S}}(k) = \left\{
\begin{matrix}
sym^k(\sO_X\oplus\sL^{-1}) & \hbox{if} & i=0 \;\; \hbox{and} \;\; k\ge 0, \cr
0 & \mbox{if} & i=0 \;\; \hbox{and} \;\; k< 0, \cr
0 & \mbox{if} & i=1 \;\; \hbox{and} \;\; k> -2, \cr
(sym^{-k-2}(\sO_X\oplus\sL^{-1}))^{\vee}\otimes\sL & \mbox{if} & i=1 \;\; \hbox{and} \;\; k\le -2.
\end{matrix} \right.
\end{eqnarray*}
\end{lem}
\ \\ \ \\

\begin{lem} \label{tangent-relative} With the same notation as above, if $\dim\; X >1$, then one has
\begin{eqnarray*}
H^1(\bar{S}, \tau_{\bar{S}/X}) =0
\end{eqnarray*}
\end{lem}

\pf Consider the short exact sequence
\begin{eqnarray} \label{euler}
0\lra \sO_{\bar{S}}\lra \pi^*(\sO_{\bar{S}}\oplus\sL)\otimes\sO(1)\lra\tau_{\bar{S}/X}\lra0
\end{eqnarray}

This sequence induces a long exact sequence of cohomology

\begin{eqnarray*}
\cdots\lra H^1(\bar{S}, \pi^*(\sO\oplus\sL)\otimes\sO(1))\lra H^1(\bar{S}, \tau_{\bar{S}/X})
\lra H^2(\bar{S}, \sO)\lra\cdots
\end{eqnarray*}
\\
But $R^1\pi_*\sO(k)=0$ for $k\ge 0$, so using the Leray spectral sequence together with the projection formula and lemma \ref{projectivebundle} one gets:
\\ \\
\begin{eqnarray*}
\begin{matrix}
H^1(\bar{S}, \pi^*(\sO\oplus\sL)\otimes\sO(1)) & = & H^1(X, (\sO\oplus\sL)\otimes R^0\bar{\pi}_*\sO(1) )\cr
 & = & H^1(X, (\sO\oplus\sL)\otimes (\sO\oplus\sL^{-1}))\cr
  & = & H^1(X, \sO\oplus\sL\oplus\sL^{-1}\oplus\sO) \hfill 
\end{matrix}
\end{eqnarray*}
and 
\begin{eqnarray*}
H^2(\bar{S},\sO_{\bar{S}}) = H^2(X, \sO_X) 
\end{eqnarray*}
and the result follows from Bott's theorem, since $\sL$ is an ample line bundle
\qed
\\ \\

\begin{lem} \label{pullback-tangent} 
\begin{eqnarray*}
H^1(\bar{S},\bar{\pi}^*\tau_X)= 0
\end{eqnarray*}
\end{lem}

\pf Again, Leray spectral sequence together with projection formula and lemma \ref{projectivebundle} give us
\begin{eqnarray*}
H^1(\bar{S},\bar{\pi}^*\tau_X) = H^1(X, \tau_X)
\end{eqnarray*}
and the results follows from Bott's theorem since the fundamental weight associated to $\tau_X$ has regularity $\ne 1$ \qed
\\ \\

\begin{lem} \label{tangent}
\begin{eqnarray*}
H^1(\bar{S}, \tau_{\bar{S}}) = 0
\end{eqnarray*}
\end{lem}
\ \\
\pf The result follows from \ref{tangent-relative} and \ref{pullback-tangent} together with the long exact sequence of cohomology associated to the short exact sequence

\begin{eqnarray} \label{relativ}
0\lra \tau_{\bar{S}/X} \lra \tau_{\bar{S}} \lra \pi^*\tau_X\lra 0
\end{eqnarray}
\hfill\qed
\\ \\

\begin{lem} \label{H2}
\begin{eqnarray*}
H^2(\bar{S}, \tau_{\bar{S}}\otimes\sO([-Z])) = 0
\end{eqnarray*}
\end{lem}
\ \\ \\
\pf By \ref{ivinskis2}, $\sO([Z]) \cong \sO(N)\otimes\bar{\pi}^*\sL^N$, therefore
\\
\begin{eqnarray*}
H^2(\bar{S}, \tau_{\bar{S}}\otimes\sO([-Z])) = H^2(\bar{S}, \tau_{\bar{S}}\otimes \sO(-N)\otimes
\bar{\pi}^*\sL^{-N})
\end{eqnarray*}
\ \\ \ \\
Tensoring the short exact sequence (\ref{relativ}) with $\sO(-N)\otimes \bar{\pi}^*\sL^{-N}$ one sees that the vanishing of $H^2(\bar{S}, \tau_{\bar{S}}\otimes\sO([-Z]))$ follows from that of 
$
H^2(\bar{S}, \tau_{\bar{S}/X} \otimes \sO(-N)\otimes \bar{\pi}^*\sL^{-N}) 
$
and
$
H^2(\bar{S}, \bar{\pi}^*(\tau_X) \otimes \sO(-N)\otimes \bar{\pi}^*\sL^{-N}) 
$.
\\ \\
Projection formula together with \ref{projectivebundle} says that \\
\begin{eqnarray*}
\begin{matrix}
H^2(\bar{S}, \bar{\pi}^*(\tau_X) \otimes \sO(-N)\otimes \bar{\pi}^*\sL^{-N}) & = &
H^1(X, \tau_X \otimes \sL^{-N}\otimes R^1\bar{\pi}_*\sO(-N)) \hfill \cr
 & = & H^1(X, \tau_X \otimes \sL^{-N}\otimes (sym^{N-2}(\sO_X\oplus\sL^{-1}))^{\vee}\otimes\sL)
\end{matrix}
\end{eqnarray*}
\\
Now the vanishing of this group follows from the vanishing of the groups \\
\begin{eqnarray*}
H^1(X, \tau_X \otimes \sL^j) \hskip1cm \hbox{for}\;\; -N< j <-1;
\end{eqnarray*}
\\
which in turn are zero because of Bott's theorem.
\\ \\
A similar argument (short exact sequence \ref{euler} together with projection formula and lemma \ref{projectivebundle}) shows that 
$
H^2(\bar{S}, \tau_{\bar{S}/X} \otimes \sO(-N)\otimes \bar{\pi}^*\sL^{-N})  = 0
$. \hfill\qed
\\ \\

\section{ Kuranishi space of deformations of cyclic coverings}

The infinitesimal Torelli theorem for cyclic coverings could be true for trivial reasons, i.e., if the corresponding Kuranishi space is trivial or discrete. In this section we will show that it is not the case in general. \label{kuranishi} In doing so, we closely follows an argument by Wavrik \cite{wavrik}, who showed the non triviality of the Kuranishi space of deformations for general cyclic coverings of the proyective space.
\\ \\
Let $X=\z$ with $\p$ a maximal parabolic subgroup, $Z$ be a simple covering of $X$ as before and $B$ be the Kuranishi space of deformations of $Z$. Then $Z\subset \bar{S}=\pP(\sO_X\oplus\sL^{-1})$ is a simple covering of $X$ as in the previous section and we have  short exac sequences

\begin{eqnarray} \label{canonical}
o\lra \sO_{\bar{S}}(-Z)\lra \sO_{\bar{S}} \lra i_*\sO_Z \lra 0
\end{eqnarray}

\begin{eqnarray} \label{restriction}
0\lra \tau_Z \lra \tau_{\bar{S}}\otimes\sO_Z \lra {\sN} \lra 0
\end{eqnarray}
\ \\
and

\begin{eqnarray} \label{relativ}
o\lra \tau_{\bar{S}/X} \lra \tau_{\bar{S}} \lra \pi^*\tau_X\lra 0
\end{eqnarray}
\ \\
where $\tau$ stands for the tangent bundle and $\sN$ stands for the normal bundle. 
\\ \\
Equation \ref{canonical} gives rice, after tensoring with $\tau_{\bar{S}}$, to the short exac sequence

\begin{eqnarray} \label{tensor}
0 \lra \tau_{\bar{S}}\otimes \sO([-Z]) \lra \tau_{\bar{S}} \lra i_*(\tau_{\bar{S}}\otimes\sO_Z) \lra 0
\end{eqnarray}
\\ \\
Associated to \ref{restriction} and \ref{tensor}  there are long exact sequences

\begin{eqnarray} \label{long restriction}
\cdots\lra H^0(Z, \sN) \lra H^1(Z, \tau_Z) \lra H^1(Z, \tau_{\bar{S}}\otimes\sO_Z)\lra \cdots 
\end{eqnarray}

\begin{eqnarray} \label{long relativ}
{\small \cdots\rightarrow H^1(\bar{S}, \tau_{\bar{S}})\rightarrow H^1(\bar{S}, i_*(\tau_{\bar{S}}\otimes\sO_Z)) \rightarrow
H^2(\bar{S}, \tau_{\bar{S}}\otimes \sO([-Z]))\rightarrow\cdots}
\end{eqnarray}
\\ \\
We would like to have $H^1(Z, \tau_{\bar{S}}\otimes\sO_Z)=0$ in \ref{long restriction}, in order to compute the dimension of
$H^1(Z, \tau_Z)$. Since
$H^1(Z, \tau_{\bar{S}}\otimes\sO_Z)=H^1(\bar{S}, i_*(\tau_{\bar{S}}\otimes\sO_Z))$, by \ref{long relativ}, we need to show that
\begin{eqnarray*}
H^1(\bar{S}, \tau_{\bar{S}})=H^2(\bar{S}, \tau_{\bar{S}}\otimes \sO([-Z]))=0
\end{eqnarray*}
which is the content of \ref{tangent} and \ref{H2}, so $H^1(Z, \tau_{\bar{S}}\otimes\sO_Z)=0$.
\\ \\
Therefore we have a short exact sequence
\begin{eqnarray*}
H^0(Z, \tau_{\bar{S}}\otimes\sO_Z)\lra H^0(Z, \sN)\lra H^1(Z, \tau_Z)\lra 0
\end{eqnarray*}
and $\dim H^1(Z, \tau_Z) >0$ provided the first map is not surjective. The following propositions tell us that this will be the case in general and are due to Wavrik (see \cite{wavrik}), though he stated them only for $\z = \pP^n $, but the proofs are the same, {\it mutatis mutandis}. \\ \\
\begin{prop}
\begin{eqnarray*}
h^0(Z,\sN) = {\displaystyle \sum_{j=0}^N h^0(\z, \sL^j) -1}.
\end{eqnarray*}
\end{prop} \ \\ 
\begin{prop}
$H^0(Z, \tau_{\bar{S}}\otimes \sO_Z) = 0$ if $d(N-1)> d_0$, where $d=c_1(\sL)$ and $d_0$ is the constant appearing in remark \ref{canonico}, which for ${\bf G}$ of type $A_n$ is equal to $n+1$.
\end{prop} \ \\
 In conclusion, in the general situation the Kuranishi space of deformations of cyclic coverings of $\z$ is not discrete.
\ \\
\section{ Infinitesimal Torelli for cyclic coverings of $\z$.}

Let $X=\z$ with $\p$ a maximal parabolic subgroup, $Z$ be a simple covering of $X$ as before and $B$ be the Kuranishi space of deformations of $Z$. Suppose B is smooth 
and let $m:B\to D$ be the local period map into the period domain $D$
of all Hodge structures on $H^n(Z,\C)$, where $n=\dim Z$. The map $m$ is
holomorphic, with tangent map at the point $Z$:
\begin{eqnarray*}
T:H^1(Z,\tau_Z)\to T_{m(0)}D=\oplus_p 
\Hom(H^{n-p}(Z,\Dif^p_Z),H^{n-p+1}(Z,\Dif^{p-1}_Z))
\end{eqnarray*}
We say that the infinitesimal Torelli theorem holds for $Z$ if the tangent
map $T$ is injective.
\ \\
Obviously it will be enought to show that there is a $p$ such that the map
\begin{eqnarray*}
T_p:H^1(Z,\tau_Z)\to \Hom(H^{p}(Z,\Dif^{n-p}_Z),H^{p+1}(Z,\Dif^{n-p-1}_Z))
\end{eqnarray*}
obtained from $T$ and the natural projection into the $p$-factor, is injective.
In particular, it will suffy to show it for $p=0$ or, what amounts to the same
because of Serre's duality, it will be enough to show that the map
\begin{eqnarray}\label{surjection}
H^0(Z, \dif_Z )\otimes H^{n-1}(Z,\Dif^1_Z )\lra H^{n-1}(Z,\Dif^1_Z \otimes\dif_Z )
\end{eqnarray}
is surjective, where as usual, $\dif_Y=\Dif^n_Y$
for any smooth $n$-dimensional variety $Y$.
\ \\ \\ 
For $Z$ smooth, the exact sequence
\begin{eqnarray}\label{hartshorne}
0\lra\sO_Z(-Z)\lra\Dif^1_{\bar{S}|Z}\stackrel{\beta}{\lra}\Dif^1_Z\lra 0
\end{eqnarray}
induces exact sequences
\begin{eqnarray}
0\lra\Dif^{n-p}_Z\otimes\sO_Z(-Z)\stackrel{\alpha}{\lra}\Dif^{n-p+1}_{\bar{S}|Z}
\stackrel{\beta}{\lra}\Dif^{n-p+1}_Z\lra 0\; ,
\end{eqnarray}
which for every locally free $\sO_Z$-module $\sG$ give rise, after tensoring 
with an appropriate multiple of $\sO_Z(Z)$, 
to a complex $(K^{\bullet},d_{\bullet})$, where
\begin{eqnarray}
K^j:=\Dif^{n-p+j+1}_{\bar{S}|Z}\otimes\sO_{Z}((j+1)Z)\otimes\sG
\end{eqnarray}
and the derivatives are induced by $\alpha\circ\beta$. This 
complex is
exact, except at $j=0$, where $\ker\; d_0=\Dif^{n-p}_Z\otimes\sG$, therefore 
the complex 
$(K^{\bullet},d_{\bullet})$ is quasi-isomorphic to $\Dif^{n-p}_Z\otimes\sG$.\\ 
Moreover, the exact sequence \ref{hartshorne} implies
$K^{p}=\dif_{Z}\otimes\sO_{Z}(pZ)\otimes\sG$ and $K^j=0$ for 
$j> p$.\\ \\
Observe that $K^{p}[-p]=F^{p}K^{\bullet}$, so that we have a map
$\H^p(Z,K^p[-p])\stackrel{\gamma}{\lra} \H^{p}(Z,K^{\bullet})$.
\\ \\
Following \cite{ivinskis}, let us define the map $\gamma_p(\sG)$ as
the composition
\begin{eqnarray*}
H^0(Z,\dif_{Z}\otimes\sO_{Z}(pZ)\otimes\sG)\cong 
\H^p(K^p[-p])\stackrel{\gamma}{\lra} \H^{p}(K^{\bullet})\cong
H^p(Z,\Dif^{n-p}_Z\otimes\sG).
\end{eqnarray*}
The following lemma gives us sufficient conditions for the map \ref{surjection}
to be surjective.
\begin{lem}[\cite{ivinskis}]\label{ivinskis3}
If the map
\begin{eqnarray*}
H^0(Z,\dif^2_Z\otimes\sO_Z((n-1)Z)) \stackrel{\gamma_{n-1}(\dif_Z)}{\lra} 
H^{n-1}(Z,\Dif^1_Z \otimes\dif_Z )
\end{eqnarray*}
and the map
\begin{eqnarray*}
H^0(Z, \dif_Z )\otimes H^0(Z,\dif_Z\otimes\sO_Z((n-1)Z)) 
\lra H^0(Z, \dif_Z )\otimes H^{n-1}(Z,\Dif^1_Z )
\end{eqnarray*}
are surjective, then the map \ref{surjection} is surjective and the 
infinitesimal Torelli theorem holds for $Z$.
\end{lem}

\pf
One has a commutative diagram
\begin{eqnarray}
\begin{matrix}
H^0(Z,\dif_Z)\otimes H^{n-1}(Z,\Dif^1_Z) & \lra & H^{n-1}(Z,\Dif^1_Z\otimes\dif_z) \cr
 & & \cr
{\tiny id\otimes\gamma_{n-1}(\sO_Z)} \uparrow & & \uparrow 
{\tiny id\otimes\gamma_{n-1}(\dif_Z)} \cr
& & \cr
H^0(Z,\dif_Z)\otimes H^0(Z,\dif_z\otimes\sO_Z((n-1)Z) & \stackrel{\eta}{\lra} & 
H^0(Z,\dif^2_Z\otimes\sO_Z((n-1)Z))
\end{matrix}
\end{eqnarray}
and  the assumptions of the  lemma imply the surjectivity of the map in the top row.
\ \\ \ \\
\begin{thm} \label{torelli}
Let $X=\z$, $Z\lra X$ be a simple covering of degree $N$ with respect to $\sO_X(d)$ and $d_0,\mu$ 
be as in remark \ref{canonico} and lemma \ref{bott1}. If
$d(N-1)-d_0> \mu$ or $d(N-1)-d_0>n-1$, then the infinitesimal Torelli theorem holds for $Z$.
\end{thm}

\pf We will freely use  the notation of section \ref{coverings}, with
$\sL=\sO_X(d)$. Remember that $n=dim\; X$.\\ \\
According with \ref{ivinskis3}, it  is  enough to show the surjectivity of the maps
\begin{eqnarray*}
H^0(Z,\dif^2_Z\otimes\sO_Z((n-1)Z)) \stackrel{\gamma_{n-1}(\dif_Z)}{\lra} 
H^{n-1}(Z,\Dif^1_Z \otimes\dif_Z )
\end{eqnarray*}
and 
\begin{eqnarray*}
H^0(Z, \dif_Z )\otimes H^0(Z,\dif_Z\otimes\sO_Z((n-1)Z)) 
\lra H^0(Z, \dif_Z )\otimes H^{n-1}(Z,\Dif^1_Z )
\end{eqnarray*}
\ \\ \ \\
I) The surjectivity of $\gamma_p(\dif_Z)$ follows from the 
spectral sequence $H^i(Z,K^j)\Lra \sH^{i+j}(K^{\bullet})$ as soon as 
$H^{s+1}(Z,K^{p-s-1})=0$ for $s=0, \cdots , p-1$. In our case $p=n-1$,
$K^j:=\Dif^{n-p+j+1}_{\bar{S}|Z}\otimes\sO_{Z}((j+1)Z)\otimes\sG$ and $\sG=\dif_Z$,
in other words, we need to show that
\begin{eqnarray*}
H^{s+1}(Z, \Dif^{n-s}_{\bar{S}|z}\otimes\sO_Z((n-1-s)Z)\otimes\dif_Z)=0
\end{eqnarray*}
for $s=0, \cdots , n-2$. \\
The classical short exact sequence for the relative differentials
\begin{eqnarray*}
0\lra \bar{\pi}^*\Dif^1_X\lra \Dif^1_{\bar{S}}\lra \Dif^1_{\bar{S}/X}\lra 0
\end{eqnarray*}
induces, after exterior product and restriction to $Z$, the short exact sequence
\begin{eqnarray*}
0\lra f^*\Dif^{n-s}_X\lra \Dif^{n-s}_{\bar{S}|Z}\lra f^*(\Dif^{n-s-1}_X\otimes \sL^{-1})
\lra 0
\end{eqnarray*}
for every $s$, which in turn produces, after tensoring with 
$\sO_Z((n-1-s)Z)\otimes\dif_Z$, the exact sequence
\begin{eqnarray*}
\begin{matrix}
0& \lra & f^*(\Dif^{n-s}_X\otimes \sL^{(n-s)\cdot N-1}\otimes\dif_X ) & \lra & 
\Dif^{n-s}_{\bar{S}|Z}\otimes \sO_Z((n-1-s)Z)\otimes\dif_Z & & \cr
 &  & & & & & \cr
 &      &         & \lra &
f^*(\Dif^{n-s-1}_X\otimes \sL^{(n-s)\cdot N-2}\otimes\dif_X ) & \lra & 0.
\end{matrix}
\end{eqnarray*}
Projection formula together with lemma \ref{ivinskis1} gives
\begin{eqnarray*}
H^{s+1}(Z, f^*(\Dif^{n-s}_X\otimes \sL^{(n-s)\cdot N-1}\otimes\dif_X)) =
\oplus_{i=0}^{N-1}H^{s+1}(X,\Dif^{n-s}_X\otimes \sL^{(n-s)\cdot N-1-i}\otimes\dif_X)
\end{eqnarray*}
and also
\begin{eqnarray*}
H^{s+1}(Z,f^*(\Dif^{n-s-1}_X\otimes \sL^{(n-s)\cdot N-2}\otimes\dif_X )) =
\oplus_{i=0}^{N-1}H^{s+1}(X,\Dif^{n-s-1}_X\otimes \sL^{(n-s)\cdot N-2-i}\otimes\dif_X )
\end{eqnarray*}
\ \\
Observe that $n-s\ge 2$ and $0\le i\le N-1$, so that $(n-s)N-1-i\ge (n-s)N-2-i \ge N-1$. 
Since $\sL^j=\sO_X(jd)$ for every $j$, then $d[(n-s)N-2-i]\ge d(N-1)$. Finally, as observed in
remark \ref{canonico}, $\omega_X=\sO_X(-d_0)$,
therefore the assumption on $d$ and $N$, together 
with proposition \ref{bott1}
give us the vanishing of this cohomology groups and therefore the vanishing of
\begin{eqnarray*}
H^{s+1}(Z, \Dif^{n-s}_{\bar{S}|z}\otimes\sO_Z((n-1-s)Z)\otimes\dif_Z)
\end{eqnarray*}
which implies the surjectivity of $\gamma_{n-1}(\dif_Z)$, as desired. \\ \\
II) Let us now consider the map
\begin{eqnarray*}
H^0(Z, \dif_Z )\otimes H^0(Z,\dif_Z\otimes\sO_Z((n-1)Z)) 
\lra H^0(Z, \dif_Z )\otimes H^{n-1}(Z,\Dif^1_Z )
\end{eqnarray*}
The surjectivity of this map will follows from the surjectivity of the map
\begin{eqnarray*}
H^0(Z,\dif_Z\otimes\sO_Z((n-1)Z)) 
\lra H^{n-1}(Z,\Dif^1_Z ).
\end{eqnarray*}
Now, for every $p$ we have short exact sequences
\begin{eqnarray*}
0 \lra \Omega^{n-p-1}_{Z}(-(k+1)Z)\lra \Omega^{n-p}_{S|Z}(-kZ)\lra \Omega^{n-p}_{Z}(-kZ)\lra 0 ,
\end{eqnarray*}
which induce long exact sequences
{\small
\begin{eqnarray*}
\begin{matrix}
\cdots   \lra  H^p\left(\Omega^{n-p}_{Z}(-(n-p-1)Z)\right)  & 
\begin{matrix}
\beta_p \cr \lra \cr \;
\end{matrix}  
  & H^{p+1}\left(\Omega^{n-p-1}_{Z}(-(n-p-2)Z)\right)\lra \cr
 & \lra & H^{p+1}\left(\Omega^{n-p}_{S|Z}(-(n-p-1)Z)\right)\lra \cdots \;\; ,
 \end{matrix}
\end{eqnarray*}
}
so the surjectivity of $\beta_p$ for every $p$ will give us the surjectivity we are looking for, but this will be the
case if $H^{p+1}\left(\Omega^{n-p}_{S|Z}(-(n-p-1)Z)\right)=0$ for every $p$, which follows from the computations on I), since
$d_0\ge 1$ always.
\\ \\

\section*{Acknowledgments}

This work started while the first named author was at the University 
Duisburg-Essen in a sabatical year. He wants to thank 
H\'el\`ene Esnault for her invitation, as well as
all other members of the group in Essen, for their hospitality and the 
marvelous ambient they have built there, from the human as well as from the
mathematical point of view. We are also in debt with Vicente Navarro, who asked about the non triviality of the Kuranishi space  of deformations of cyclic coverings, leading to section \ref{kuranishi}. Finally we are deeply in debt to Eckart Viehweg, who suggested the problem for Grassmannians long ago. 

\bibliographystyle{plain}
\renewcommand\refname{References}

{}
\end{document}